\documentclass{amsart}
\usepackage{amssymb,amsmath,amsthm,epsf,epsfig,dsfont,bbm}

\usepackage{latexsym}
\usepackage{upref, eucal}
\usepackage[all]{xy}

\newcommand {\nc} {\newcommand}
\newcommand {\enm} {\ensuremath}

\def \d{\delta}

\nc {\bdm} {\begin{displaymath}}
\nc {\edm} {\end{displaymath}}

\newtheorem {theorem} {\bf{Theorem}}[section]
\newtheorem {lemma}[theorem] {\bf Lemma}
\newtheorem {proposition}[theorem] {\bf Proposition}

\newtheorem {corollary}[theorem] {\bf Corollary}
\numberwithin {equation}{section}

\newcommand{\Ou}{\enm{\mathcal{O}}}

\newcommand{\A}{\enm{\mathbb{A}}}

\nc{\J}{\enm{\mathcal{J} }}
\nc {\Z} {\enm{\mathbb{Z}}}
\nc {\form}[1] {\enm{\mbox{\underline{for}}}_{#1}}
\nc {\prol}[1] {\enm{\mbox{\underline{prol}}_{{#1}^*}}}

\nc {\stk} {\stackrel}
\newcommand{\imp}{\Rightarrow}
\newcommand{\map}{\rightarrow}

\newcommand{\beqar}{\begin{eqnarray*}}
\newcommand{\eeqar}{\end{eqnarray*}}
\newcommand{\inj}{\hookrightarrow}
\newcommand{\Pn}[2] {\ensuremath{ {\mathbb{P}}^{#1}_{#2}}}
\nc{\Quot}[3]{\enm{ {\mathfrak{Quot}_{ {#1}/{#2}/{#3}}}}}
\nc{\Hilb}[2]{\enm{ {\mathfrak{Hilb}_{ {#1}/{#2}}}}}
\newcommand{\mfrak}[1]{\mathfrak{#1}}

\newcommand{\bb}[1]{\mathbb{#1}}
\newcommand{\mcal}[1]{\mathcal{#1}}

\nc {\Coh}[4] {\ensuremath{H^{#1}(\Pn{#2}{},{#3}({#4}))}}
\nc {\Ch}[3] {\enm{H^{#1}(X_t,{#2}_t({#3}))}}
\nc {\Qphi}[4]{\enm{ {\mathfrak{Quot}^{~#4}_{ {#1}/{#2}/{#3}}}}}
\nc {\Gra}[4]{\enm{ {\mathfrak{Grass}_{#2}({#3},{#4})}}}
\nc {\HomA}[2]{\enm{\mathrm{Hom}_A{#1}{#2}}}
\nc {\tr}{\mathrm{tr}}

\nc {\C}[2]{\enm{\left(\begin{array}{l} {#1} \\ {#2} \end{array} \right)}}
\nc {\mat}[4]{\enm{\left(\begin{array}{ll}{#1} & {#2} \\ {#3} & {#4}
\end{array}\right)}}

\def \mb{\mbox}

 \def \Z{{\mathbb Z}}

   \def \h{\hat{\ }}

\def \d{\delta}

\def \R1{R((q))[q']\h}

\DeclareMathOperator{\Spec}{\mathrm{Spec}}

\newcommand{\Hom}{\mathrm{Hom}}

\title{On Frobenius and Fibers of Arithmetic Jet Spaces}
\author{James Borger \and Arnab Saha}
\date{}

\begin{document}
\maketitle

\begin {abstract}
In this article, given a scheme $X$ we show the existence of canonical 
lifts of Frobenius maps in 
an inverse system of schemes obtained from the fiber product of the canonical 
prolongation sequence of arithmetic jet spaces
$J^*X$ and a prolongation sequence $S^*$ over the scheme $X$. 
As a consequence, for any smooth group scheme $E$, if $N^n$ denote the 
kernel of the canonical projection map of the $n$-th jet space $J^nE \map E$, 
then the inverse system $\{N^n\}_n$ is a prolongation sequence. 
\end {abstract}

\section{Introduction}
The purpose of this short note is to make an observation which is a
generalisation of a result shown in \cite{drin} for Drinfeld modules. 
Let $B$ be a Dedekind domain and fix a maximal ideal $\mfrak{p} \in \Spec B$
with $k:=B/\mfrak{p}$ a finite field and let $q=|k|$. 
Let $R$ be the $\mfrak{p}$-adic completion of $B$. Denote by $\mfrak{m}$ the 
maximal ideal of the complete, local ring $R$ and $\iota: B \inj R$ the natural
inclusion. Then let $\pi \in B$ be such that $\iota(\pi)$ generates the maximal 
ideal $\mfrak{m}$ in $R$. Since $\iota$ is an injection, by 
abuse of notation, we will consider $\pi$ as an element of $B$ as well.
Then $k\simeq R/(\pi)$. Do note here that the identity 
map on $R$ lifts the $q$-power Frobenius on $R/(\pi)$. Let $V=\Spec R$.

To motivate the main result in
this article, let $E$ be a group scheme over $V$. Then by \cite{bor11b,
bui95}, one can consider, for all integer $n \geq 0$, the $n$-th jet 
space $J^nE$ (Here by $J^nE$, we understand the algebraic jet space in
\cite{bor11b}. Buium's jet spaces are $p$-adic formal schemes obtained by 
taking 
formal completions of $J^nE$). Due to functorial reasons, $J^nE$ is also
a group scheme for all $n$. There are canonical maps, the projection map
 $u:J^nE \map J^{n-1}E$ and the lift of Frobenius with respect to $u$ denoted
$\phi:J^nE \map J^{n-1}E$ for all $n\geq 1$. The map $\phi$ is 
associated to a $\pi$-derivation $\d$ on the structure sheaves which will
be defined in section $2$ in detail. Such a system of schemes is called a 
{\it prolongation sequence}.
Let $J^*E$ denote the prolongation sequence
of jet spaces.
Then for each $n$, we have the following short exact sequence of group 
schemes
\begin{align}
\label{short1}
0 \longrightarrow N^n \longrightarrow J^nE \stk{u}{\longrightarrow} E 
\longrightarrow 0
\end{align}
where $N^n$ is the kernel of $u$. Then it is easy to see that the projection
map $u:J^nE \map J^{n-1}E$ induces a $u:N^n \map N^{n-1}$ by restriction. 
Let $N^*$ denote the inverse system of such schemes.
However, the $\d$ or $\phi$ does not restrict to $N^*$  because that 
would imply that $E$ has a lift of Frobenius which is false in general.

Let $X$ be a scheme over $V$. 
For any prolongation sequences $T^*$ and $S^*$ over $V$ with morphisms 
$T^0 \map X$ and $S^0 \map X$, let us define $T^* \times_X S^* :=\{T^n \times_X 
S^n\}_{n=0}^\infty$. Note that $T^* \times_X S^*$ is not  apriori a 
prolongation sequence but is an inverse system with the projection maps induced
from the projection maps from $T^*$ and $S^*$. For any inverse system of 
schemes $U^*$ with maps
$u:U^n \map U^{n-1}$ for all $n \geq 1$, we will say $U^*$ admits a lift
of Frobenius if for all $n$, there are maps $\phi:U^n \map U^{n-1}$ which 
are lifts of Frobenius with respect to the projection maps $u$.
Our main result is
\begin{theorem}
\label{final}
If $X$ is a scheme over $V$, $J^*X \times_X S^*$ admits a canonical lift of 
Frobenius. 
\end{theorem}

We call this canonical lift of Frobenius as {\it Lateral Frobenius} and
is denoted by $\mfrak{f}:J^nX \times_X S^n \map J^{n-1}X\times_X S^{n-1}$
for all $n \geq 1$.
Then if $E$ is a smooth group scheme, then $J^nE$ are also smooth and 
hence flat. Therefore a lift of Frobenius is equivalent to a $\pi$-derivation
in this case.
It is now easy to see that by choosing $X=E$ and $S^*$ to be the constant 
prolongation sequence given by $S^n:=V$ for all $n$, with the given map $S^0
\map E$ as the identity section of the group scheme $E$, it follows that 
$N^n$ admits a lift of Frobenius and thus is a prolongation sequence. 

\section{Prelimineries}
\subsection{Witt vectors}
Here we recall some basic facts about Witt vectors and arithmetic jet spaces.
Witt vectors over a general Dedekind domain with finite residue fields
 was developed in \cite{bor11a}. 
For the sake of our article, we will briefly review the general construction.
Let $B$ be a Dedekind domain and fix a maximal ideal $\mfrak{p} \in \Spec B$
with $k:=B/\mfrak{p}$ a finite field and let $q=|k|$. 
Let $R$ be the $\mfrak{p}$-adic completion of $B$. Denote by $\mfrak{m}$ the 
maximal ideal of the complete, local ring $R$ and $\iota: B \inj R$ the natural
inclusion. Then let $\pi \in B$ be such that $\iota(\pi)$ generates the maximal 
ideal $\mfrak{m}$ in $R$. Since $\iota$ is an injection, by 
abuse of notation, we will consider $\pi$ as an element of $B$ as well.
Then $k\simeq R/(\pi)$. 

Do note here that the identity 
map on $R$ lifts the $q$-power Frobenius on $R/(\pi)$. We will now review
the theory of $\pi$-typical Witt vectors over $R$ with maximal ideal 
$\mfrak{m}$. All the rings in this section are $R$-algebras.

Let $C$ be an $A$-algebra with structure map $u:A \map C$. In this paper,
any ring homomorphism $\psi: A \map C$ will be called the {\it lift of 
Frobenius} if it satisfies the following:

(1) The reduction mod $\pi$ of $\psi$ is the $q$-power Frobenius, that is,
$\psi(x) \equiv u(x)^q$ mod $ \pi C$.

(2) The restriction of $\psi$ to $R$ is identity.

Let $C$ be an $A$-algebra with structure map $u:A \map C$. A $\pi$-derivation
$\d$ from $A$ to $C$ means a set theoretic map satisfying the following for
all $x,y \in A$

\beqar
\label{der}
\d(x+y) &=& \d (x) + \d (y) +C_\pi(u(x),u(y)) \\
\d(xy) &=& u(x)^q \d (y) +  u(y)^q \d (x) + \pi \d (x) \d (y) \\
\eeqar
such that $\d$ when restricted to $R$ is $\d(r) = (r-r^q)/\pi$ for
all $r \in R$ and 
$$C_\pi(X,Y) = 
 \left\{\begin{array}{ll} 0, & \mb{ if $R$ is positive characteristic}\\
\frac{X^q+Y^q - (X+Y)^q}{\pi}, & \mb{ otherwise}
\end{array} \right. $$
It follows that the map $\phi: A \map C$ defined as 
$$
\phi(x) := u(x)^q + \pi \d (x)
$$ 
is an $R$-algebra homomorphism and is a lift of the Frobenius.
Considering this operator $\d$ leads 
to Buium's theory of arithmetic jet spaces \cite{bui95,bui00,
bui09}.

Note that this definition depends on the choice of uniformizer $\pi$, but in a transparent way:
if $\pi'$ is another uniformizer, then $\d(x)\pi/\pi'$ is a $\pi'$-derivation, and this correspondence
induces a bijection between $\pi$-derivations and $\pi'$-derivations.

Given an $R$-algebra $A$, the ring of $\pi$-typical Witt vectors 
$W(A)$ can be defined as
the unique $R$-algebra $W(A)$ with a $\pi$-derivtion $\d$ on 
$W(A)$ such that, given any $R$-algebra $C$ with a 
$\pi$-derivation $\d$ on it and an $R$-algebra map $f:C \map
A$, there exists an unique $R$-algebra homomorphism $g:C \map W(A)$ satisfying-
$$\xymatrix{
W(B) \ar[d] & \\
A & C \ar[l]_f \ar[ul]_g 
}$$
and $g$ satisfies $g \circ \d = \d \circ g$.
In \cite{bor11a} (following the approach of \cite{joyal} to the usual 
$p$-typical Witt vectors),
the existence of such a $W(A)$ is shown and that it is also 
obtained from the classical definition of Witt vectors using ghost vectors.

\subsection{Prolongation sequences and Jet spaces}
Let $V= \Spec R$ and $X$ and $Y$ be schemes over $V$. We say a
pair $(u,\d)$ is a {\it prolongation}, and write 
$Y \stk{(u,\d)}{\map} X$, if $u: Y \map X$ is a map of 
schemes over $V$ and $\d: \Ou_X \map u_*\Ou_Y$ is a $\pi$-derivation 
satisfying 
$$\xymatrix{
B \ar[r] &  u_* \Ou_Y \\
B \ar[u]^\d \ar[r] &  \Ou_X \ar[u]_\d \\
} $$

Following \cite{bui00}, a {\it prolongation sequence} is a sequence of prolongations
$$
\xymatrix{
V & T^0 \ar_-{(u,\d)}[l] & T^1 \ar_-{(u,\d)}[l] & \cdots\ar_-{(u,\d)}[l]},
$$
where each $T^n$ is a scheme over $V$. We will often use the notation
$T^*$ or $\{T_n\}_{n \geq 0}$.
Note that if the  $T^n$s are flat over $V$ then a 
$\pi$-derivation $\d$ is equivalent to a lift of Frobenius $\phi$ as defined 
above.

Prolongation sequences form a category $\mcal{C}_{V^*}$, where a morphism $f:T^*\to U^*$ is 
a family of morphisms $f^n:T^n\to U^n$ commuting with both the $u$ and $\d$, in the evident sense.
This category has a final object $V^*$ given by $V^n=\Spec R$ for all $n$, 
where each $u$ is the identity and
each $\d$ is the given $\pi$-derivation on $R$.

For any $V$-scheme $Y$,  for all $n \geq 0$ we define 
the $n$-th jet space $J^nX$ (relative to $V$) as 
$$
J^nX (Y) := \Hom_\d(W_n^*(Y),X)
$$
where $W_n^*(Y)$ is defined as in \cite{bor11b}. We will not define $\Hom_\d
(W_n^*(Y),X)$ in full generality here. Instead, we will define 
$\Hom_\d(W_n^*(Y),X)$
in the affine case which is obviously simpler but will offer an
intuitive understanding of the definition.
Let our schemes be affine and $X = \Spec A$ and $Y=\Spec C$. Then 
$W_n^*(Y)= \Spec W_n(C)$ and $\Hom_\d(W_{n}^*Y,X)= \Hom_\d(A, W_n(C))$, where
$W_n(C)$ is the ring of truncated Witt vectors of length $n+1$ and 
$f \in \Hom_\d(A,W_n(C))$ is a ring homomorphism that satisfies the following
$$\xymatrix{
A \ar[r]^-f & W_n(C) \\
R \ar[r]_-{\exp_\d} \ar[u] & W_n(R) \ar[u]
}$$
where $R\stk{\exp_\d}{\longrightarrow} W_n(R)$ is the universal map from the 
definition of Witt vectors.

Then $J^*X:= \{J^nX \}_{n \geq 0}$ forms a prolongation sequence and is 
called the {\it canonical prolongation sequence}. By \cite{bui00},
 $J^*X$ satisfies the following 
universal property---for any $T^* \in \mcal{C}_{V^*}$ and $X$ a
 scheme, we have 
$$\Hom(S^0,X) \simeq \Hom_{\mcal{C}_{S^*}}(S^*, J^*X)$$

\section{ Construction of the Lateral Frobenius}
Let all our schemes be over a base $V$ which has a lift of Frobenius. 
Let $S^*= \{S^n\}_{n=0}^\infty$ be a 
prolongation sequence with a morphism $S^0 \map X$. Also, given a 
prolongation sequence, let $S^{*-1} := \{S^{n-1}\}_{n=1}^\infty$. 
For any prolongation sequence $T^*$ with a morphism $T^0 \map X$, let us define 
$T^* \times_X S^* :=\{T^n \times_X 
S^n\}_{n=0}^\infty$. Note that $T^* \times_X S^*$ is not  apriori a 
prolongation sequence.

\subsection{Affine $N$-space case.}
Let $X = \A^N$ where $N$ may represent an arbitrary cardinality. 
For any $V$-scheme $T$, let $\prod_j^{j+m} T:= \underbrace{T \times_V
\cdots \times_V T}_{m+1-\mb{times}}$.
Then 
by definition of jet spaces, we have a canonical isomorphism of schemes 
$J^nX \simeq W_n^N$, where $W_n$ is the truncated Witt vectors
of length $n+1$. Therefore as a scheme, $J^nX$ can be 
canonically identified as $J^nX \simeq \prod_0^n X$.
Hence $J^nX \times_X S^n \simeq (\prod_1^nX)\times_V S^n$.
Let $w: J^nX \map \prod_0^n X$ be the ghost map of the product of Witt vectors 
and the right-hand side is referred to as the ghost components.

We define the {\it lateral Frobenius} $\mfrak{f}$ as the unique morphism which
makes the following diagram commutative
$$\xymatrix{
J^nX \times_X S^n \ar[d]_{\mfrak{f}} \ar[r]^{w \times \mathbbm{1}} & 
(\prod_0^n X) \times_X S^n \ar[r]^-{\sim} &
(\prod_1^n X) \times_V S^n \ar[d]^{\mfrak{f}} \\
J^{n-1}X \times_X S^{n-1} \ar[r]^{w \times \mathbbm{1}} & 
(\prod_0^{n-1} X) \times_X S^{n-1} \ar[r]^-{\sim}
& (\prod_1^{n-1} X) \times_V S^{n-1} 
}$$
where $\mfrak{f}: (\prod_1^n X) \times_V S^n \map (\prod_1^{n-1}X) \times_V 
S^{n-1}$
is given by the left-shift operator on the ghost components
 $\mfrak{f}((w_1,\cdots, w_n),s) = ((w_2,...,w_n),\phi(s))$ and completely
determines the 
map $\mfrak{f}:J^nX \times_X S^n \map J^{n-1}X \times_X S^{n-1}$. The map 
$\mfrak{f}$ 
is unique and is given by $\mfrak{f}(x_1, \cdots, x_n,s)= (z_1,\cdots ,z_{n-1},
\phi(s))$
where $(x_1,\cdots ,x_n) \mapsto (z_1,\cdots, z_{n-1})$ is the Frobenius map 
of Witt vectors and hence the map is a lift of Frobenius as well. It is clear 
that $\mfrak{f}$ behaves functorially with respect to $X$.

Let $l: J^nX \times_X S^n \map J^nX \times_V S^n$ 
denote the natural map induced from the structure map $X \map V$ for all $n$. 

\begin{proposition}
The lateral Frobenius $\mfrak{f}$ satisfies the following commutative diagram
for all $n \geq 2$
$$\xymatrix{
J^nX \times_X S^n \ar[d]_-{\mfrak{f}} \ar[r]^l & J^nX \times_V S^n 
\ar[d]^-{\phi \times \phi} \\
J^{n-1}X \times_X S^{n-1} \ar[dr]_-{(\phi \times \phi) \circ l}
 & J^{n-1}X \times_V S^{n-1} \ar[d]^-{\phi \times \phi} \\
& J^{n-2}X \times_V S^{n-2} 
}$$
\end{proposition}
{\it Proof.} This follows by checking commutativity on the ghost component
since $w$ is injective$\qed$

\subsection{The case of a general affine scheme}
Let $S^*= \{S^n\}_n$ be a prolongation sequence with a morphism 
$S^0\stk{a}{\map} X$. We will denote this data as $S^* \stk{a}{\map} X$. Let 
$X$ be an affine scheme. Then it satisfies an equaliser diagram 
$$\xymatrix{
X \ar[r]^-f & Y \ar@<1ex>[r]^-g \ar[r]_-h & Z \\
}$$
where $Y=\bb{A}^N$ and $Z=\bb{A}^M$, $N$ and $M$ are arbitrary and $f$ is 
injective.
Since the jet space functor $J^n$ preserves equalisers \cite{bor11b} we obtain
an equaliser  of the corresponding jet spaces with the following commutative 
diagram 

$$\xymatrix{
& S^* \ar[dl]_a \ar@{=}[rr] & & S^* \ar@{=}[rr] \ar[dl]_{f\circ a} & & S^* 
\ar[dl]_{g\circ f \circ a} \\
X \ar[rr]^-f & &  Y \ar@<1ex>[rr]^-g \ar[rr]_-h & & Z & \\
J^*X \ar[rr]^-f \ar[u]^u & & J^*Y \ar[u]^u \ar@<1ex>[rr]^g \ar[rr]_-h & & 
J^*Z \ar[u]^u & \\
}$$

\begin{lemma}
The map
$J^*X \times_Y S^* \stk{f}{\map} J^*Y \times_Y S^*$ is injective.
\end{lemma}
{\it Proof.} If $f(x,s)=f(x',s') \iff (f(x), s) = (f(x'),s')$ satisfying
$u \circ f(x) = f\circ a (s) = u \circ f(x') = f \circ a(s')$ in $Y$. This
implies that $x=x'$ since $f$ is injective and $s=s'$ and we are done. $\qed$

\begin{lemma}
We have $J^*X \times_Y S^* = J^*X \times_X S^*$
\end{lemma}
{\it Proof.} We have $(x,s) \in J^*X \times_Y S^* \iff f \circ \pi(x) = f \circ
a(s)$ in $Y \iff \pi(x) = a(s)$ in $X~ (\text{Since } f$ is injective$) 
\iff (x,s) \in J^*X \times_X S^*$ and we are done. $\qed$

Combining the above two results we get an injective map $J^*X \times_X S^* 
\stk{f}{\map} J^*Y \times_Y S^*$. 

\begin{proposition}
The following 
$$\xymatrix{
J^*X \times_X S^* \ar[r]^-f & J^*Y \times_Y S^* \ar[r]<1ex>^g \ar[r]_h 
& J^*Z \times_Z S^*}$$ is an equaliser diagram.
\end{proposition}
{\it Proof.} We need to show exactnes, that is, if for a $(y,s) \in J^*Y 
\times_Y S^*$ such that $g(y,s) = h(y,s)$ then we claim that 
$(y,s) \in J^*X \times_X S^*$. Now $g(y,s) = h(y,s) \iff (g(y),s) = (h(y),s)
\iff g(y) = h(y)$ in $J^*Z \iff y \in J^*X$ (Since 
$J^*$ preserves an equaliser diagram). Therefore $(y,s) \in J^*X \times_X S^*$
and we are done. $\qed$

\begin{corollary}
\label{affine}
If $X$ is affine, then $J^*X \times_X S^*$ admits the lift of Frobenius 
induced from $\mfrak{f}$.
\end{corollary}
{\it Proof.} Since $Y = \bb{A}^N$ and $Z=\bb{A}^M$, we have shown that
$J^*Y\times_Y S^*$ and $J^*Z \times_Z S^*$ admit the lift of Frobenius
$\mfrak{f}$. Then the result follows from the universal property of equalisers
$$\xymatrix{
J^nX \times_X S^n \ar@{.>}[d]^{\mfrak{f}} \ar[r] & J^nY \times_Y S^n 
\ar[d]^{\mfrak{f}} \ar@<1ex>[r]^f \ar[r] & J^nZ \times_Z S^n \ar[d]^{\mfrak{f}}
 \\
J^{n-1}X \times_X S^{n-1} \ar[r] & J^{n-1}Y \times_Y S^{n-1}
\ar@<1ex>[r]^f \ar[r] & J^{n-1}Z \times_Z S^{n-1} \\
}$$
$\qed$

\subsection{The case of a general scheme}

Let $X$ be a scheme which admits a finite cover by affines $X_i$. Let 
$Y = \coprod X_i$ and $Z = Y \times_X Y$. 
A prolongation sequence such that $S^*\stk{a}{\map} X$, induces the map 
$S^* \map J^*X$ of prolongation sequences by the universal property of 
jet spaces. Define a new prolongation sequence $\tilde{S}^*$ given by
$\tilde{S}^n := S^n \times_{J^nX} J^nY$ for all $n$. Define 
$\bar{S}^*:= \tilde{S}^* \times_{S^*} \tilde{S}^*$. Then we have the following
commutative diagram of coequalisers.

$$\xymatrix{
& \bar{S}^* \ar[dl]_a \ar@<1ex>[rr]^-g \ar[rr]_-h& &  \tilde{S}^* 
\ar[dl]_{g\circ a} \ar[rr]^-f &  & S^* \ar[dl]_{f \circ g \circ a} & \\
Z \ar@<1ex>[rr]^-g \ar[rr]_-h & & Y \ar[rr]^-f & & X &\\
J^*Z \ar@<1ex>[rr]^-g \ar[u]^u \ar[rr]_-h & & J^*Y \ar[u]^u \ar[rr]^-f & & 
J^*X \ar[u]^u &\\
}$$

\begin{proposition}
\label{etale}
$J^*Y \times_Y \tilde{S}^* \map J^*X \times_X S^*$ is an etale surjection, 
that is, $J^nY \times_Y \tilde{S}^n \map J^nX \times_X S^n$ is an 
etale surjection for all $n$.
\end{proposition}
{\it Proof.} By \cite{bor11b}, since $Y \map X$ is etale, $J^*Y \map J^*X$ and
$\tilde{S}^* \map S^*$ are etale as well.
Therefore, $J^*Y \times_Y \tilde{S}^* \map J^*X \times_X S^*$ is etale.

For any $V$-scheme $T$, if $Q \map V$ is the generic point, then denote
$T_Q:= T \times_V Q$. If $P \map V$ is the closed point, then denote 
$\overline{T}:= T \times_V P$.
Now for all $n$, $(J^nY \times_Y \tilde{S}^n)_Q = (\prod_1^n Y_Q) \times 
\tilde{S}^n_Q$ which clearly surjects over $(J^nX \times S^n)_Q= (\prod_1^n X_Q)
\times S^n_Q$. By \cite{bui95}, we have $\overline{J^nY} = \overline{J^nX} 
\times_{\overline{X}} \overline{Y}$. Therefore,
$\overline{J^nY} \times_{\overline{Y}} \overline{\tilde{S}^n}= \overline{J^nX}
 \times_{\overline{X}} \overline{Y} \times_{\overline{Y}} \overline{\tilde{S}^n}
 = \overline{J^nX} \times_{\overline{X}} \overline{\tilde{S}^n}$ which 
clearly surjects over $\overline{J^nX} \times_{\overline{X}} \overline{S^n}$. 
Hence 
$J^nY \times_Y \tilde{S}^n \map J^nX \times_X S^n$ is a surjection of schemes
for all $n$ and we are done. $\qed$

\begin{proposition}
\label{fprod}
$(J^*Y \times_Y \tilde{S}^*) \times_{J^*X \times_X S^*} 
(J^*Y \times_Y \tilde{S}^*) = (J^*Y \times_{J^*X} J^*Y) \times_{(Y \times_X Y)}
(\tilde{S}^* \times_{S^*} \tilde{S}^*)$
\end{proposition}
{\it Proof.} Let $b= g \circ a$. 
Define a map $(J^*Y \times_Y \tilde{S}^*) \times_{J^*X \times_X 
S^*} (J^*Y \times_Y \tilde{S}^*) \map (J^*Y \times_{J^*X} J^*Y) 
\times_{(Y \times_X Y)} (\tilde{S}^* \times_{S^*} \tilde{S}^*)$ as 
$((y,\tilde{s}),(y',\tilde{s}')) \mapsto ((y,y'),(\tilde{s},\tilde{s}'))$. 
We claim that the above map is well-defined. Consider $((y,\tilde{s}),
(y',\tilde{s}')) \in (J^*Y \times_Y \tilde{S}^*) \times_{J^*X \times_X S^*} 
(J^*Y \times_Y \tilde{S}^*)$. Then it satisfies the following-
$(i)~ u(y)= b(\tilde{s}), u(y')= b(\tilde{s'})$ in $Y$ and $(ii)~ 
f(y,\tilde{s})=f(y',\tilde{s}') \imp f(y)=f(y')$ in $J^*X$ and $f\circ 
b(\tilde{s}) = f \circ b(\tilde{s}')$. But then the above two conditions
implies $(i)~f(y) = f(y')$ in $J^*X$ and $f(\tilde{s})=f(\tilde{s'})$ in 
$S^*$ and $(ii)~ u(y) = b(\tilde{s}),~ u(y') = b(\tilde{s}')$ which implies 
$((y,y'),
(\tilde{s},\tilde{s'})) \in (J^*Y \times_{J^*X} J^*Y) \times_{(Y \times_X Y)}
(\tilde{S}^* \times_{S^*} \tilde{S}^*)$ and hence proves the claim of 
well-definedness. Clearly this map has an inverse and we are done. $\qed$

\begin{corollary}
$$\xymatrix{
J^*Z \times_Z \bar{S}^*
\ar[r]<1ex>^-g \ar[r]_-h
& J^*Y \times_Y \tilde{S}^* \ar[r]^-f & J^*X \times_X S^*
}$$ 
is a coequaliser diagram.
\end{corollary}
{\it Proof.} Since $\bar{S}^* = \tilde{S}^* \times_{S^*} \tilde{S}^*$, the
result follows from proposition \ref{etale} and \ref{fprod}. $\qed$


{\it Proof of Theorem \ref{final}}
Since $Z$ and $Y$ are affine, by corollary \ref{affine}, we know that 
$J^*Z\times_Z S^*$ and $J^*Y \times_Y S^*$ admit the lift of Frobenius 
$\mfrak{f}$ and hence the result for $J^*X \times_X S^*$ follows from the 
universal property of coequalisers,
$$\xymatrix{
J^nZ \times_Z \bar{S}^n \ar[d]^{\mfrak{f}} \ar[r]_h \ar@<1ex>[r]^g & 
J^nY \times_Y \tilde{S}^n \ar[d]^{\mfrak{f}} \ar[r]^-f& J^nX \times_X S^n 
\ar@{.>}[d]^{\mfrak{f}} \\
J^{n-1}Z \times_Z \bar{S}^{n-1} \ar[r]_h \ar@<1ex>[r]^g
 & J^{n-1}Y \times_Y \tilde{S}^{n-1}
\ar[r]^f & J^{n-1}X \times_X S^{n-1} \\
}$$
$\qed$


\begin{thebibliography}{1}

\bibitem{bor11a}
J.~Borger.
\newblock The basic geometry of Witt vectors I: The affine case.
\newblock {\em Algebra \& Number Theory}, 5(2):231--285, 2011.

\bibitem{bor11b}
J.~Borger.
\newblock The basic geometry of Witt vectors II: Spaces.
\newblock {\em Mathematische Annalen}, 351(4):877--933, 2011.

\bibitem{drin}
J.~Borger and A.~Saha.
\newblock Differential Characters of Drinfeld Modules and de Rham Cohomology.
\newblock {\it http://arxiv.org/abs/1703.05677.}

\bibitem{bui95}
A.~Buium.
\newblock Differential characters of abelian varieties over $p$-adic fields.
\newblock {\em Inventiones mathematicae}, 122(1):309--340, 1995.

\bibitem{bui00}
A.~Buium.
\newblock Differential modular forms.
\newblock {\em Journal fur die Reine und Angewandte Mathematik}, 520:95--168,
  2000.

\bibitem{bui09}
A.~Buium and B.~Poonen.
\newblock Independence of points on elliptic curves arising from special points
  on modular and shimura curves, II: local results.
\newblock {\em Compositio Mathematica}, 145(03):566--602, 2009.

\bibitem{joyal}
Andr{\'e} Joyal.
\newblock $\delta$-anneaux et vecteurs de Witt.
\newblock {\em CR Math. Rep. Acad. Sci. Canada}, 7(3):177--182, 1985.

\end{thebibliography}

\end{document}